\documentclass{conm-p-l}

\copyrightinfo{2002}{American Mathematical Society}

\newtheorem{theorem}{Theorem}[section]
\newtheorem{lemma}[theorem]{Lemma}
\newtheorem{prop}[theorem]{Proposition}
\newtheorem{cor}[theorem]{Corollary}

\theoremstyle{definition}
\newtheorem{definition}[theorem]{Definition}
\newtheorem{example}[theorem]{Example}

\theoremstyle{remark}
\newtheorem{remark}[theorem]{Remark}

\numberwithin{equation}{section}

\begin{document}

\title[Symplectic Supermanifolds and Courant algebroids]{On the structure of graded
symplectic supermanifolds and Courant algebroids.}

\author{Dmitry Roytenberg}
\address{Department of Mathematics, Pennsylvania State University, University Park, PA
16802}
\email{roytend@math.psu.edu}

\subjclass[2000]{}
\date{February 26, 2002}
\subjclass{Primary 53D05, 81T70; Secondary 51P05, 81T45}
\keywords{supermanifolds, symplectic, algebroid, BRST}

\begin{abstract}
This paper is devoted to a study of geometric structures expressible in terms
of graded symplectic supermanifolds. We extend the classical BRST formalism
to arbitrary pseudo-Euclidean vector bundles \( E\rightarrow M_{0} \) by canonically
associating to such a bundle a graded symplectic supermanifold \( (M,\Omega ) \),
with \( \textrm{deg}(\Omega )=2 \). Conversely, every such manifold arises
in this way. We describe the algebra of functions on \( M \) in terms of \( E \)
and show that ``BRST charges'' on \( M \) correspond to Courant algebroid
structures on \( E \), thereby constructing the standard complex for the latter
as a generalization of the classical BRST complex. As an application of these
ideas, we prove the acyclicity of ``higher de Rham complexes'', a generalization
of a classic result of Fr\"{o}hlicher-Nijenhuis, and derive several easy but
useful corollaries. 
\end{abstract}
\maketitle

\section{Introduction.}

\newcommand{\bform}{<\cdot ,\cdot >}

\newcommand{\A}{{\mathcal{A}}}

Graded symplectic supermanifolds have been known to physicists since the 70s,
providing framework for the so-called BRST formalism. Such supermanifolds were
obtained by adjoining odd generators (``ghosts'' and ``antighosts'') to
the classical algebra of observables, and having them satisfy canonical Poisson
bracket relations. In the simplest case, such supermanifolds are of the form
\( M=\Pi ({\mathfrak {g}}\oplus {\mathfrak {g}}^{*})\times T^{*}M_{0} \), where
\( M_{0} \) is the configuration space acted upon by a Lie algebra \( {\mathfrak {g}} \),
and \( \Pi  \) denotes parity shift. The symplectic form is \( \Omega =dp_{i}dq^{i}+d\theta _{a}d\xi ^{a} \),
where the coordinates \( (\xi ^{a}) \) on \( \Pi {\mathfrak {g}} \) are the
``ghosts'', those \( (\theta _{a}) \) on \( \Pi {\mathfrak {g}}^{*} \) --
the ``antighosts''. One assigns the ghost degree \( +1 \) to the ghost variables,
while the antighosts are assigned degree \( -1 \). The symplectic form \( \Omega  \)
has total ghost degree zero, so that one can replace \( T^{*}M_{0} \) by an
arbitrary symplectic manifold \( P \) acted upon by \( {\mathfrak {g}} \).
The classical BRST algebra is the graded Poisson algebra \( \wedge ^{\cdot }({\mathfrak {g}}\oplus {\mathfrak {g}}^{*})\otimes C^{\infty }(P) \)
of functions on \( M \). One then constructs an odd self-commuting Hamiltonian
of ghost degree \( +1 \) (the ``BRST charge'') which acts on the BRST algebra
via the Poisson bracket as a differential; its cohomology (the BRST cohomology)
is identified with the ``physical'' observables.

A detailed mathematical study of the above special case of classical BRST algebra
and its quantization was undertaken by Kostant and Sternberg \cite{KoSt} who
related the BRST cohomology with Lie algebra cohomology. They showed, in particular,
that if the action of \( {\mathfrak {g}} \) on \( P \) is ``nice'', the
BRST cohomology is concentrated in degree \( 0 \) and is isomorphic to the
algebra of functions on the reduced phase space.

The classical BRST formalism can be extended and generalized in various directions.
For instance, to treat more complicated symmetries it is required to further
extend the BRST algebra by including ``ghosts for ghosts'' etc. The homological
algebra becomes quite involved in this case (see \cite{Kje}). The purpose of
this paper is to generalize the BRST formalism to the case of an arbitrary pseudo-Euclidean
vector bundle \( E \) over a manifold \( M_{0} \), and to study the resulting
geometric structures, Courant algebroids. Consider the special case \( M_{0}=\textrm{point} \),
treated in \cite{KoSt}. Here we have a pseudo-Euclidean vector space \( (V,\bform ) \),
\( M=\Pi V \), \( \Omega =\frac{1}{2}d\xi ^{a}g_{ab}d\xi ^{b} \), where \( g_{ab}=<e_{a},e_{b}> \)
for some basis \( (e_{a}) \) of \( V \). Not having a ghost-antighost splitting,
we assign degree \( +1 \) to all \( \xi ^{a} \)'s, so that the degree of \( \Omega  \)
is now \( +2 \). The exterior algebra \( \wedge ^{\cdot }V^{*} \) becomes
a graded Poisson algebra; the Poisson bracket has degree \( -2 \). Hence, the
space \( \wedge ^{2}V^{*} \) of quadratic Hamiltonians becomes a Lie algebra,
isomorphic to \( {\mathfrak {so}}(V) \). It acts on all Grassman polynomials
\( \wedge ^{k}V^{*} \) via the Poisson bracket. The situation is entirely analogous
to the well-known realization of the Lie algebra \( {\mathfrak {sp}}(V) \)
of a symplectic vector space \( (V,\omega ) \) as the lie algebra of quadratic
Hamiltonians \( S^{2}V^{*} \) under the canonical Poisson bracket.

If \( V \) carries a Lie algebra structure such that \( \bform  \) is ad-invariant,
one forms the \emph{Cartan structure tensor} \( \Theta \in \wedge ^{3}V^{*} \)
by setting \( \Theta (X,Y,Z)=<[X,Y],Z> \). It plays the role of a ``BRST charge''
for the standard complex of the Lie algebra \( V \): it satisfies the structure
equation \( \{\Theta ,\Theta \}=0 \) (equivalent to the Jacobi identity), and
\( \delta =\{\Theta ,\cdot \} \) is the Chevalley-Eilenberg differential. The
Lie bracket on \( V \) is recovered as a so-called \emph{derived bracket}:
\( [X,Y]=\{\{X,\Theta \},Y\} \) (we have used \( \bform  \) to identify \( V \)
with \( V^{*} \)). This formalism was exploited by Lecomte-Roger \cite{LecRog}
and Kosmann-Schwarzbach \cite{KS3} who used it to study the homological algebra
of Lie bialgebras and quasi-Lie bialgebras, respectively.

If \( M_{0}\times V \) is a trivial pseudo-Euclidean vector bundle over \( M_{0} \),
one sets \( M=\Pi V\times T^{*}M_{0} \), with \( \Omega =dp_{i}dq^{i}+\frac{1}{2}d\xi ^{a}g_{ab}d\xi ^{b} \).
One needs to assign degree \( 2 \) to the momenta \( p_{i} \) so that \( \Omega  \)
is homogeneous of degree \( 2 \). Thus it is proper to denote \( M=V[1]\times T^{*}[2]M_{0} \),
to indicate the shift in grading. Note that with this choice of grading every
homogeneous function on \( M \) has a nonnegative degree compatible with its
parity: \( M \) is a so-called \emph{N-manifold} \cite{Sev2}. It turns out
that one can canonically associate a symplectic N-manifold \( M \) of degree
\( 2 \) to any, not necessarily trivial, pseudo-Euclidean vector bundle \( E\rightarrow M_{0} \),
and that conversely, every such manifold arisies in this way (Theorem \ref{thm:deg2}).
Local trivializations of \( E \) (as a pseudo-Euclidean bundle) give rise to
affine Darboux charts of the form \( V[1]\times T^{*}[2]U \), and \( M \)
is glued together from these by affine coordinate transformations (\ref{eqn:coord_change})
corresponding to changes of trivialization of \( E \). An observation, often
overlooked, is that all canonical constructions in the category of vector bundles
should be equivariant with respect to the group of bundle automorphisms covering
diffeomorphisms of the base (as opposed to just those that leave every point
of the base fixed). The Lie algebra of this group is known as the \emph{Atiyah
algebra}, consisting of so-called covariant differential operators (CDO's) on
\( E \). In case of pseudo-Euclidean vector bundles, one restricts to the subalgebra
of operators preserving \( \bform  \). This algebra is realized as the Lie
algebra \( \A ^{2} \) of quadratic Hamiltonians on \( M \), under the Poisson
bracket. It acts on all polynomial functions on \( M \), preserving the degree.

Whereas in studying structures on vector spaces (such as Lie algebras, bialgebras,
etc.) the use of supermanifold formalism is entirely optional and a matter of
one's taste, it offers significant advantages already when we go to vector bundles.
The analogue of a Lie algebra in the world of vector bundles is a \emph{Lie
algebroid}. A Lie algebroid structure on a vector bundle \( A\rightarrow M_{0} \)
is usually described in terms of an anchor map \( A\rightarrow TM_{0} \) and
a Lie bracket of sections, but it can be equivalently described as a derivation
\( d_{A} \) of \( \Gamma (\wedge ^{\cdot }A^{*}) \) of degree \( +1 \) and
square zero. This differential can be viewed as a homological vector field on
the N-manifold \( A[1]=\Pi A \) (in the terminology of \cite{Sev2}, an \emph{NQ-structure}).
The standard example is the tangent bundle \( TM_{0} \) with identity anchor
and the commutator bracket of vector fields; the complex in this case is just
the de Rham complex of \( M_{0} \). This homological description facilitates
an elegant and concise formulation of such important notions as modules and
homomoprhisms of Lie algebroids. The classical definition of homomorphism is
rather complicated, especially when different base manifolds are involved; the
classical notion of a module, while certainly reasonable, fails to include perhaps
the most important module of all, the adjoint module (see \cite{Mac} for a
thorough discussion). The homological formulation overcomes these difficulties
\cite{Vain}. Its potential is yet to be fully realized.

As soon as a pseudo-Euclidean structure on a vector bundle \( E\rightarrow M_{0} \)
is introduced, one quickly realizes that any reasonable notion of ``ad-invariance''
is incompatible with the notion of Lie algebroid, except in the trivial case
of a bundle of Lie algebras. One has to sacrifice the skew-symmetry of the bracket
of sections by adding a symmetric part which is, in a suitable sense, infinitesimal
(a coboundary), or else introduce anomalies into the Jacobi and other basic
identities. Thus \emph{Courant algebroids} are born. The first basic example
(called \emph{standard}), on \( E=TM_{0}\oplus T^{*}M_{0} \), is due to T.
Courant \cite{Cou}, who used it to describe integrability of so-called \emph{Dirac
structures} on \( M_{0} \) (including 2-forms, bivectors and distributions)
in a unified way. The general definition appeared in \cite{LWX1} as a result
of a ``doubling'' construction for \emph{Lie bialgebroids} \cite{MacXu}.
This leads, for a Lie bialgebroid \( (A,A^{*}) \), to a Courant algebroid structure
on \( E=A\oplus A^{*} \), generalizing Courant's example. \v{S}evera \cite{Sev}
showed how Courant algebroids appear naturally in the context of two-dimensional
variational problems. He also classified so-called \emph{exact Courant algebroids},
deformations of the standard one by closed 3-forms. These were later used to
describe twisted Poisson manifolds with a closed 3-form background (\cite{SevWe},
see also \cite{Roy3-QuasiLie}).

To describe Courant algebroids homologically, the (merely Poisson, not symplectic)
N-manifold \( E[1]=\Pi E \) is not enough; one needs its minimal symplectic
realization \( (M,\Omega ) \) mentioned above. It turns out (Theorem \ref{thm:CA})
that Courant algebroid structures on \( E \) are in 1-1 correspondence with
``BRST charges'' on \( M \), i.e. cubic Hamiltonians \( \Theta  \) satisfying
the structure equation \( \{\Theta ,\Theta \}=0 \). The BRST differential \( D=\{\Theta ,\cdot \} \),
acting on the algebra \( \A ^{\cdot } \) of polynomial functions on \( M \),
gives the standard complex for the Courant algebroid. For Lie bialgebroids this
description (the ``homological double'') was given in \cite{Roy1} (in that
case \( M=T^{*}\Pi A \) with an appropriate assignment of grading). In case
\( M_{0} \) is a point, we recover Cartan's structure tensor. At this point
we should remark that instead of N-manifolds we could consider general graded
manifolds \cite{Vor3}; all our constructions extend to this case at the small
extra cost of having to consider the even and odd symplectic forms separately,
as the degree is no longer tied to parity. This way we would obtain even and
odd Courant superalgebroids over supermanifolds. We stick with N-manifolds for
simplicity of presentation.

One remarkable feature of Courant algebroids, apparent in the homological formulation,
is that they can be used as target spaces for a general class of 3-dimensional
topological field theories. Special cases include the Chern-Simons gauge theory
(when \( M_{0} \) is a point) and Park's topological membrane theory \cite{Park}
(when \( E \) is an exact Courant algebroid). This fits in the general Batalin-Vilkovisky
framework described in \cite{AKSZ}. The applications to the topology of
3-manifolds are beginning to be investigated; the program was outlined to
us by Park \cite{ParkPrivate}, who also considers symplectic NQ-manifolds of
higher degree (see \cite{Park} for application to deformation quantization).

When \( E=TM_{0}\oplus T^{*}M_{0} \) is the standard Courant algebroid, \( M=T^{*}\Pi TM_{0} \),
and \( \Theta =\Theta _{0} \) is just the Hamiltonian corresponding to the
de Rham vector field \( d \) on \( \Pi TM_{0} \). The standard (BRST) complex
gets decomposed into a series of subcomplexes \( (\A ^{k,\cdot },D) \) which
we call the \emph{higher de Rham complexes}. For \( k=0 \) this is just the
de Rham complex of \( M_{0} \), whereas \( \A ^{1,\cdot }[1] \) can be identified
with vector fields on \( \Pi TM_{0} \), with \( D=[d,\cdot ] \); in general,
\( \A ^{k,l} \) consists of the symbols of differential operators on \( \Pi TM_{0} \)
of order \( k \) and degree \( l-k \); \( (\A ^{k,\cdot },D) \) can be interpreted
as the standard complex of the Lie algebroid \( TM_{0} \) with coefficients
in the \( k \)th exterior power of the adjoint module.

We prove (Proposition \ref{prop:acyclicity}) that the higher de Rham complexes
are acyclic for \( k>0 \) by constructing a chain homotopy for \( D \) as
a vector field \( \iota  \) on \( M \) of degree \( -1 \). This generalizes
the classic result of Fr\"{o}licher and Nijenhuis \cite{FroNij} on the structure
of derivations of the ring of differential forms. Several corollaries, potentially
of independent interest, follow immediately from this: Cor. \ref{cor:rigidityTM}
(the rigidity of the standard Lie algebroid structure on \( TM_{0} \)), Cor.
\ref{cor:Poisson} (any Lie algebroid structure on \( T^{*}M_{0} \) compatible
with the standard one on \( TM_{0} \) comes from a Poisson structure on \( M_{0} \))
and Cor. \ref{cor:standardCA_cohomology} (the cohomology of the standard Courant
algebroid is isomorphic to the de Rham cohomology of \( M_{0} \)). This last
corollary implies, in particular, that \v{S}evera's exact Courant algebroids
are the only nontrivial deformations of the standard one.

This paper is organized as follows. In Section \ref{sec:GradedMan} we review
the generalities on graded manifolds and N-manifolds, following \cite{Vor3}
and \cite{Sev2}; in Section \ref{sec:deg1and2} we study the structure of symplectic
N-manifolds of degree 1 and 2 in detail; Section \ref{sec:CA} is devoted to
the homological (BRST) description of Courant algebroids, their standard complexes
and cohomology; finally, in Section \ref{sec:exact} we concentrate on the special
case of the standard Courant algebroid and its deformations.

We would like to thank Anton Alekseev, Jae-Suk Park, Pavol \v{S}evera, Ted Voronov,
Alan Weinstein and Ping Xu for useful discussions and advice. We also thank
the Erwin Schr{\"o}dinger Institute where some of this work was carried out,
for hospitality.

\section{\label{sec:GradedMan} Graded Manifolds and N-manifolds.}

A \emph{graded manifold} is a (super)manifold with an additional grading in
the structure sheaf. Specifically, it is a manifold \( M \) that possesses
a coordinate atlas (of so-called \emph{affine charts}) in which each local coordinate
is assigned a \emph{degree} (or \emph{weight}), and the coordinate transformations
are required to preserve the total weight. In this note we shall require all
weights to be integers, although it is not a necessary restriction, in principle.
The grading can also be conveniently described by means of the \emph{Euler vector
field} 
\[
\epsilon =\sum w(x^{i})x^{i}\frac{\partial }{\partial x^{i}},\]
 where \( w(x^{i}) \) denotes the weight of \( x^{i} \). The weights are just
the eigenvalues with respect to the action of \( \epsilon  \): we say that
a function \( f \) is \emph{homogeneous of weight} \( k \) if \( \epsilon \cdot f=kf \).
We shall denote the space of all such functions by \( \A ^{k} \) and call the
graded algebra \( \A ^{\cdot }=\bigoplus _{k\geq 0}\A ^{k} \) the \emph{algebra
of polynomial functions} on \( M \); the algebra of all smooth functions is
a completion of \( \A ^{\cdot } \). The \emph{degree of \( M \)} is by definition
the highest weight of a local coordinate.

Any (super)manifold \( M \) becomes a graded manifold if we assign zero weight
to all the coordinates (i.e. set \( \epsilon =0 \)). The total space \( M_{1} \)
of any vector bundle \( A \) over a base (super)manifold \( M_{0} \) becomes
a graded manifold if the base coordinates are assigned zero weight, while all
the fibre coordinates are assigned the same weight \( n \). The standard choice,
of course, is \( n=1 \). In fact, any manifold of degree \( 1 \) is easily
seen to be a vector bundle. 

In general, different local coordinates can have different weights, and the
grading can be thought of as a generalization of a vector bundle structure;
a comparison with physical units of measurement is also helpful. A general graded
manifold \( M \) has the following structure. Let \( \A _{k} \) denote the
(graded) subalgebra of \( \A ^{\cdot } \) locally generated by functions of
weight \( \leq k \). Then \( \A ^{\cdot } \) is filtered by the \( \A _{k} \)'s:
\[
\A _{0}\subset \A _{1}\subset \A _{2}\subset \cdots \subset \A _{d}=\A ^{\cdot },\]
 where \( d \) is the degree of \( M \). Clearly, \( \A ^{k}=\A _{k}/\A _{k+1} \).
Correspondingly, we have a tower of fibrations 
\[
M_{0}\leftarrow M_{1}\leftarrow M_{2}\leftarrow \cdots \leftarrow M_{d}=M\]
 Here \( M_{0} \) is a (super)manifold, \( M_{1} \) is a vector bundle over
\( M_{0} \), and for \( k\geq 1 \), \( M_{k}\leftarrow M_{k+1} \) is an affine
fibration. Each \( M_{k} \) is itself a graded manifold (of degree \( k \)),
and their Euler vector fields are related by the projections. \( M_{0} \) also
embeds into each \( M_{k} \) as the set of zeros of the Euler vector field;
it corresponds to the ideal \( \A ^{\geq 1} \).

In particular, \( \A _{0}=\A ^{0}=C^{\infty }(M_{0}) \), and all \( \A _{k} \)'s
and \( \A ^{k} \)'s are modules over \( \A _{0} \). It is best to think of
these as (locally free) sheaves of \( C^{\infty }(M_{0}) \)-modules over \( M_{0} \),
as they can be restricted to arbitrary open subsets of \( M_{0} \). Thus, \( \A ^{1} \)
is the sheaf of sections of some vector bundle over \( M_{0} \), while \( \A _{1}=S^{\cdot }\A ^{1} \),
etc. In what follows we shall frequently commit an abuse of notation by using
\( \A ^{k} \) to refer to both the vector bundle and its sheaf of sections.

A (non-negatively integer) graded supermanifold \( M \) is an \emph{N-manifold}
if the integer grading is compatible with parity (the underlying \( {\mathbb {Z}}_{2} \)-grading
in the structure sheaf). In plain English, this means that even coordinates
have even weights, while odd ones have odd weights. The parity operator is related
to the Euler vector field via \( \Pi =(-1)^{\epsilon } \). Thus we have an
action of the multiplicative semigroup \( ({\mathbb {R}},\times ) \) on \( M \):
\( -1 \) acts as \( \Pi  \), \( 0 \) picks out \( \A ^{0} \), while \( \lambda >0 \)
acts as \( \lambda ^{\epsilon } \). Then \( M_{0}=0\cdot M \) is an ordinary
manifold, while \( M_{1}=\Pi A \) for some ordinary vector bundle over \( M_{0} \)
(i.e. \( \A _{1}=\Gamma (\wedge ^{\cdot }A^{*}) \)). Conversely, given a vector
bundle \( A\rightarrow M_{0} \), we can define the N-manifold \( A[n] \) by
assigning degree \( n \) to the fibre coordinates; it is equal to either \( A \)
or \( \Pi A \), depending on the parity of \( n \). The \emph{support} of
an N-manifold \( M \) is the ordinary graded submanifold \( M_{\textrm{even}} \)
corresponding to the ideal \( \A ^{\textrm{odd}} \).

\begin{remark}
\label{rem:OnGrading} More generally, if \( N \) is a graded manifold and
\( E\rightarrow N \) is a graded vector bundle, we denote by \( E[n] \) the
graded manifold obtained by shifting the fibre degrees by \( n \). In particular,
\( TN \) and \( T^{*}N \) are graded vector bundles, where the Euler vector
field acts as a Lie derivative. With respect to this induced grading, the ``velocities''
have the same degree as the corresponding coordinates on \( N \), whereas the
``momenta'' have the opposite degrees. If \( N \) is an N-manifold, we may
also have to shift the parity in the fibers so that \( E[n] \) is again an
n-manifold. For example, if \( A\rightarrow M_{0} \) is an ordinary vector
bundle, \( T[1]A[1] \) is \( \Pi T(\Pi A) \) with coordinates \( (x^{i},\xi ^{a},dx^{i},d\xi ^{a}) \)
of weights \( w(x^{i})=0 \), \( w(\xi ^{a})=1 \), \( w(dx^{i})=1 \), \( w(d\xi ^{a})=2 \),
and \( T^{*}[2]A[1] \) is \( T^{*}\Pi A \) with coordinates \( (x^{i},\xi ^{a},p_{i},\theta _{a}) \)
of weights \( w(x^{i})=0 \), \( w(\xi ^{a})=1 \), \( w(p_{i})=2 \), \( w(\theta _{a})=1 \).
\end{remark}
Graded manifolds form a category: the morphisms are required to be \( \epsilon  \)-equivariant
(weight-preserving). N-manifolds form a full subcategory.

Given a graded manifold \( M \), the Euler vector field \( \epsilon  \) acts
on all canonical objects (tensors, jets, etc.) on \( M \) via the Lie derivative,
whereby these objects also acquire weights. For example, an \emph{NQ-manifold}
is an N-manifold endowed with an integrable (homological) vector field \( Q \)
of weight \( +1 \). That is, \( [\epsilon ,Q]=Q \) and \( [Q,Q]=2Q^{2}=0 \).
Such a \( Q \) is necessarily odd, so the integrability condition is nontrivial,
and \( (\A ^{\cdot },Q) \) becomes a cochain complex. The classic example of
an NQ-manifold is the anti-tangent bundle \( T[1]M=\Pi TM \) of an ordinary
manifold \( M \), where \( Q=d_{M} \) is the exterior derivative. In fact,
the assignment \( M\mapsto (T[1]M,d_{M}) \) is a functor that embeds ordinary
manifolds into NQ-manifolds as a full subcategory. Thus, NQ-manifolds should
be thought of as a natural enlargement of the category of smooth manifolds:
much of the calculus carries over.

More generally, NQ-manifolds of degree \( 1 \) are the same as Lie algebroids
over the base \( M_{0} \). For this reason it is natural to call NQ-manifolds
of degree \( n \) \emph{\( n \)-algebroids}. Interesting results and conjectures
concerning the integration of \( n \)-algebroids are contained in \cite{Sev2}.

Other objects of interest are homogeneous symplectic and Poisson structures.
For instance, a symplectic structure of degree \( n \) is a closed non-degenerate
2-form \( \Omega  \) such that \( {\mathcal{L}}_{\epsilon }\Omega =n\Omega  \).
The corresponding Poisson tensor \( W \) has weight \( -n \). We shall also
consider symplectic NQ-manifolds for which we have the additional condition
\( {\mathcal{L}}_{Q}\Omega =0 \).

Let us conclude this section with the following simple observations, showing
that a homogeneous symplectic form places severe restrictions on the structure
of \( M \):

\begin{lemma}
\label{lemma:exact}Let \( M \) be a graded manifold. Then 
\end{lemma}
\begin{enumerate}
\item Any homogeneous symplectic form \( \Omega  \) on \( M \) of weight \( n\geq 1 \)
is exact; in fact, 
\[
\Omega =d(\frac{1}{n}\iota _{\epsilon }\Omega )\]

\item Any homogeneous vector field \( V \) of weight \( m>-n \) preserving \( \Omega  \)
is Hamiltonian. In fact, 
\[
\iota _{V}\Omega =\pm d\left( \frac{1}{m+n}\iota _{V}\iota _{\epsilon }\Omega \right) \]

\end{enumerate}
\begin{proof}
Both statements are easy consequences of Cartan's homotopy formula. 
\end{proof}
\begin{example}
\label{example:T^*[1]}Given a manifold \( M_{0} \), consider \( M=T^{*}[1]M_{0}=\Pi T^{*}M_{0} \),
with \( \Omega =d\theta _{i}dx^{i}=d(\theta _{i}dx^{i}) \) of degree \( +1 \).
Functions on \( M \) correspond to multivector fields on \( M_{0} \), and
the Poisson bracket of degree \( -1 \) given by \( \Omega  \) is the Schouten
bracket. Vector fields on \( M \) of degree \( -1 \) correspond to 1-forms
on \( M_{0} \) via \( \alpha \leftrightarrow \iota _{\alpha } \). Such a vector
field preserves \( \Omega  \) if and only if \( \alpha  \) is closed; it is
Hamiltonian if and only if \( \alpha  \) is exact. This example illustates
that the second statement of the Lemma above generally fails if \( m=-n \). 
\end{example}
\begin{lemma}
\label{lemma:deg_bound}If \( (M,\Omega ) \) is a graded symplectic manifold,
the degree of \( M \) cannot exceed \( w(\Omega ) \).
\end{lemma}
\begin{proof}
In an affine chart, \( \Omega =\frac{1}{2}dx^{i}\Omega _{ij}(x)dx^{j} \), and
the total weight of each term is \( w(\Omega ) \). The statement immediately
follows from the nondegeneracy of \( \Omega  \).
\end{proof}

\section{\label{sec:deg1and2} Degree 1 and 2.}

Now let us analyze symplectic manifolds of degree 1 and 2 in detail. We shall
restrict our attention to N-manifolds, although the general case can be handled
at the small additional cost of having to consider the even and odd cases separately.

As a warm-up exercise, let us do the degree \( 1 \) case. So let \( (M,\Omega ) \)
be a symplectic N-manifold with \( w(\Omega )=1 \). According to Lemma \ref{lemma:deg_bound},
the degree of \( M \) is at most \( 1 \), so \( M=M_{1}\rightarrow M_{0} \)
is a vector bundle. More precisely, \( M=M_{1}=A^{*}[1]=\Pi A^{*} \) for some
vector bundle \( A \) over a manifold \( M_{0} \). It will be convenient to
use the Poisson bracket defined by \( \Omega  \). It is odd, nondegenerate
and has degree \( -1 \). We have \( \{\A ^{0},\A ^{0}\}=0 \), \( \{\A ^{1},\A ^{0}\}\subset \A ^{0} \)
and \( \{\A ^{1},\A ^{1}\}\subset \A ^{1} \), hence the structure is determined
by the Hamiltonian action of a Lie algebra \( \A ^{1}=\Gamma (A) \) on \( \A ^{0}=C^{\infty }(M_{0}) \).
According to the Leibniz rule, this action is given by a vector bundle map \( A\rightarrow TM_{0} \),
which is an isomorphism by the nondegeneracy of \( \{\cdot ,\cdot \} \). Finally,
by Jacobi, the Poisson bracket on \( \A ^{1} \) corresponds under this isomorphism
to the commutator of vector fields on \( M_{0} \). Thus, we have shown that
symplectic N-manifolds of degree \( 1 \) are exhausted by Example \ref{example:T^*[1]}:

\begin{prop}
\label{prop:deg1} Symplectic N-manifolds of degree \( 1 \) are in one-to-one
correspondence with ordinary smooth manifolds, via \( N\leftrightarrow (T^{*}[1]N,\Omega ) \),
where \( \Omega  \) is determined by the Schouten bracket of multivector fields.
\end{prop}
Let us notice two things. First, the Lie algebra of symmetries of \( (M,\Omega ) \)
(Hamiltonian by Lemma \ref{lemma:exact}) is \( \A ^{1} \), which is identified
with the Lie algebra of vector fields on \( M_{0} \). Second, a local coordinate
chart \( \{x^{i}\} \) on \( M_{0} \) gives rise automatically to an affine
Darboux chart \( \{x^{i},\theta _{i}\} \) on \( M \), in which \( \Omega =d\theta _{i}dx^{i}=d(\theta _{i}dx^{i}) \).

Now let us tackle the degree \( 2 \) case. By Lemma \ref{lemma:deg_bound},
such a manifold has the following structure: 
\[
M=M_{2}\rightarrow M_{1}\rightarrow M_{0}\]
 Here \( M_{0} \) is an ordinary manifold, and \( M_{1}=E[1]=\Pi E \) for
some vector bundle \( E\rightarrow M_{0} \). Let us analyze the graded Poisson
algebra of functions. The Poisson bracket is even, nondegenerate and has weight
\( -2 \). Its structure is determined by the relations \( \{\A ^{0},\A ^{0}\}=\{\A ^{0},\A ^{1}\}=0 \),
\( \{\A ^{1},\A ^{1}\}\subset \A ^{0} \), \( \{\A ^{2},\A ^{0}\}\subset \A ^{0} \),
\( \{\A ^{2},\A ^{1}\}\subset \A ^{1} \), \( \{\A ^{2},\A ^{2}\}\subset \A ^{2} \).
The first three relations, together with the symmetry and derivation properties
of the Poisson bracket, imply that \( E \) is equipped with a nondegenerate
symmetric fiberwise bilinear form \( \bform  \). In fact, \( \A _{1}=\Gamma (\wedge ^{\cdot }E^{*}) \)
is a graded Poisson subalgebra of \( \A ^{\cdot } \): the Poisson bracket is
just an extension of \( \bform  \) as a derivation in each argument. Hence,
\( M_{1} \) is a graded Poisson manifold whose symplectic leaves are the fibres
of \( \Pi E \), and the projection \( M_{2}\rightarrow M_{1} \) is a Poisson
map. One says that \( M_{2} \) is a \emph{symplectic realization} of \( M_{1} \). 

The rest of the relations define a Lie algebra structure on \( \A ^{2} \) and
its action on \( \A ^{0} \) and \( \A ^{1} \). We recall that \( \A ^{2} \)
is a locally free sheaf of \( \A ^{0} \)-modules, i.e. the sheaf of sections
of a vector bundle \( {\mathbb {A}} \) over \( M_{0} \). The Leibniz rule
implies that the action of \( \A ^{2} \) on \( \A ^{0} \) comes from a bundle
map \( a:{\mathbb {A}}\rightarrow TM_{0} \); furthermore, the Leibniz and Jacobi
identities imply that \( {\mathbb {A}} \) is actually a Lie algebroid. By nondegeneracy,
the anchor \( a \) is surjective; \( \A ^{1}\A ^{1}=\Gamma (\wedge ^{2}E^{*})\subset \A ^{2} \)
acts trivially on \( \A ^{0} \), hence the kernel of \( a \) contains \( \wedge ^{2}E^{*} \).
Furthermore, again by Leibniz and Jacobi, the action of \( \A ^{2} \) on \( \A ^{1} \)
comes from a Lie algebroid action of \( {\mathbb {A}} \) on \( E \), preserving
\( \bform  \), i.e. a Lie algebroid homomorphism from \( {\mathbb {A}} \)
into the gauge Lie algebroid of \( (E,\bform ) \); by nondegeneracy, this is
an isomorphism. Therefore, \( {\mathbb {A}} \) fits into the following exact
sequence of vector bundles over \( M_{0} \):
\begin{equation}
\label{eqn:atiyah}
0\rightarrow \wedge ^{2}E^{*}\rightarrow {\mathbb {A}}\stackrel{a}{\rightarrow }TM_{0}\rightarrow 0,
\end{equation}
 and the Lie algebroid structure on \( {\mathbb {A}} \) is that of the gauge
Lie algebroid of \\
\( (E,\bform ) \). This is the so-called \emph{Atiyah sequence} of the pseudo-Euclidean
vector bundle \( (E,\bform ) \), and the Lie algebra \( \A ^{2} \) of sections
of \( {\mathbb {A}} \) is the \emph{Atiyah algebra} of covariant differential
operators (CDO's) on \( E \) preserving \( \bform  \). The kernel \( \wedge ^{2}E^{*} \)
is a Lie algebra bundle, identified with \( {\mathfrak {so}}(E) \), the endomorphisms
of \( (E,\bform ) \) via the Poisson bracket. We conclude that the structure
of \( (M,\Omega ) \) is completely determined by \( (E,\bform ) \).

\begin{example}
\label{example:T^*[2]} (\cite{Roy1}) Let \( A\rightarrow M_{0} \) be any
vector bundle, and let \( M=T^{*}[2]A[1] \). The notation means that \( M=T^{*}\Pi A \),
with the degrees assigned in such a way that the weights of a base coordinate
and its conjugate momentum add up to \( 2 \) (Remark \ref{rem:OnGrading}).
Then the canonical symplectic form \( \Omega =dp_{i}dq^{i}+d\theta _{a}d\xi ^{a} \)
has degree \( 2 \). It is evident that in this case \( E=A\oplus A^{*} \)
, with the canonical \( \bform  \) given by \( <(X,\alpha ),(Y,\beta )>=\alpha (Y)+\beta (X) \). 
\end{example}
Conversely, suppose we are given a pseudo-Euclidean vector bundle \( (E,\bform ) \)
over a manifold \( M_{0} \). Then, of course, \( M_{1}=\Pi E=E[1] \) is a
Poisson N-manifold (of degree \( -2 \)). We shall construct its minimal symplectic
realization as follows. By Example \ref{example:T^*[2]}, \( T^{*}[2]E[1] \)
is the minimal symplectic realization of \( (E\oplus E^{*})[1] \). The map
\( E\hookrightarrow E\oplus E^{*} \) given by \( X\mapsto (X,\frac{1}{2}<X,\cdot >) \)
is an isometric bundle embedding and gives rise to a map of N-manifolds \( E[1]\hookrightarrow (E\oplus E^{*})[1] \).
We then let \( M \) be the pullback of \( T[2]E[1] \), i.e. complete the commutative
diagram 
\[
\begin{array}{ccc}
M & \longrightarrow  & T^{*}[2]E[1]\\
\downarrow  &  & \downarrow \\
E[1] & \longrightarrow  & (E\oplus E^{*})[1]
\end{array}\]
 It is clear that \( M \) is a symplectic N-submanifold of \( T^{*}[2]E[1] \),
and a minimal symplectic realization of \( M_{1}=E[1] \). Thus we have proved 

\begin{theorem}
\label{thm:deg2} Symplectic N-manifolds of degree \( 2 \) are in one-to-one
correspondence with pseudo-Euclidean vector bundles. The correspondence is given
by the above construction. Under this correspondence, the Lie algebra of degree-preserving
canonical transformations coincides with the Atiyah algebra of the bundle.
\end{theorem}
\begin{example}
If the base \( M_{0} \) is a point, then \( M=M_{1}=V[1] \) for some pseudo-Euclidean
vector space \( (V,\bform ) \). In this case the degree of \( M \) is just
\( 1 \). Given a basis \( \{e_{a}\} \) of \( V \), \( \Omega =\frac{1}{2}d\xi ^{a}g_{ab}d\xi ^{b} \),
where \( g_{ab}=<e_{a},e_{b}> \). The algebra of symmetries in this case is
\( \A ^{2}=\wedge ^{2}V^{*}\simeq {\mathfrak {so}}(V) \), acting on \( \A ^{k}=\wedge ^{k}V^{*} \),
\( k\geq 0 \), via the Poisson bracket.
\end{example}
In the general case, a choice of local coordinates \( \{x^{i}\} \) on \( M_{0} \)
and a local basis \( \{e_{a}\} \) of sections of \( E \) such that \( g_{ab}=<e_{a},e_{b}> \)
are constant gives rise to an affine Darboux chart \( (q^{i},\xi ^{a},p_{i}) \)
for \( M \): the embedding of \( M \) into \( T^{*}[2]E[1] \) is given locally
by equations \( \theta _{a}=\frac{1}{2}g_{ab}\xi ^{b} \), hence 
\[
\Omega =(dp_{i}dq^{i}+d\xi ^{a}d\theta _{a})|_{M}=dq^{i}dp_{i}+\frac{1}{2}d\xi ^{a}g_{ab}d\xi ^{b}.\]
 That is, on a local chart \( U\subset M_{0} \) such that \( E|_{U}\simeq U\times V \)
as a pseudo-Euclidean vector bundle, \( M|_{U}\simeq T^{*}[2]U\times V[1] \),
with \( \Omega  \) as above. A quadratic Hamiltonian \( H=v^{i}(q)p_{i}+\frac{1}{2}\xi ^{a}b_{ab}(q)\xi ^{b}\in \A ^{2} \)
gives rise to infinitesimal transformations 
\begin{equation}
\label{eqn:symmetry_action}
\begin{array}{ccccl}
\delta q^{i} & = & \{H,q^{i}\} & = & v^{i}\\
\delta \xi ^{c} & = & \{H,\xi ^{c}\} & = & \xi ^{a}b_{ab}(q)g^{bc}\\
\delta p_{i} & = & \{H,p_{i}\} & = & -\frac{\partial v^{j}}{\partial q^{i}}p_{j}-\frac{1}{2}\xi ^{a}\frac{\partial b_{ab}}{\partial x^{i}}\xi ^{b}
\end{array}
\end{equation}
 Integrating these, we get canonical coordinate transformations 
\begin{equation}
\label{eqn:coord_change}
\begin{array}{ccl}
q^{i} & = & q^{i}(q')\\
\xi ^{a} & = & T^{a}_{a'}(q')\xi ^{a'}\\
p_{i} & = & \frac{\partial q^{i'}}{\partial q^{i}}p_{i'}+\frac{1}{2}\xi ^{a'}\frac{\partial T^{a}_{a'}}{\partial x^{i}}g_{ab}T^{b}_{b'}\xi ^{b'}
\end{array}
\end{equation}
 corresponding to bundle transformations 
\[
\begin{array}{ccl}
x^{i'} & = & x^{i'}(x)\\
e_{a'} & = & e_{a}T_{a'}^{a}(x)
\end{array}\]
 of \( E \) such that \( T^{a}_{a'}g_{ab}T^{b}_{b'}=g_{a'b'}=\textrm{const} \).
Notice the nonlinear (affine) term in the transformation law for \( p_{i} \). 

Let us now describe the global structure of the sheaf of graded \( \A ^{0} \)-algebras
\( \A ^{\cdot } \) in terms of \( (E,\bform ) \). According to the above analysis,
it is generated by \( \A ^{1} \) in degree 1 and \( \A ^{2} \) in degree 2.
As sheaves on \( M_{0} \), \( \A ^{1}\simeq E \),\footnote{
We identify \( E \) with \( E^{*} \) via \( \bform  \)
} while \( \A ^{2}\simeq {\mathbb {A}} \). This latter is the sheaf of infinitesimal
bundle isometries of \( E \), i.e. linear vector fields on \( E \) preserving
\( \bform  \), containing \( \wedge ^{2}E\simeq {\mathfrak {so}}(E) \) as
those vector fields that leave every point of the base fixed. Therefore, 
\[
\A ^{\cdot }\simeq (\wedge ^{\cdot }E\otimes S^{\cdot }{\mathbb {A}})/I,\]
 where \( I \) is the homogeneous ideal generated by the embedding of \( \wedge ^{2}E \)
into \( {\mathbb {A}} \), i.e. by elements of the form \( \omega \otimes 1-1\otimes \omega  \),
\( \omega \in \wedge ^{2}E \). Moreover, \( \wedge ^{\cdot }E \) is a sheaf
of graded Poisson algebras, and so is \( S^{\cdot }{\mathbb {A}} \) (because
\( {\mathbb {A}} \) is a Lie algebroid), and \( I \) is a Poisson ideal by
definition of the embedding. This gives the (nondegenerate) graded Poisson algebra
structure on \( \A ^{\cdot } \).

Each homogeneous component \( \A ^{n} \) has the structure of a filtered vector
bundle over \( M_{0} \). For \( n=2k \), \( k\geq 0 \), we have 
\[
L_{k}\subset L_{k-1}\subset \cdots \subset L_{1}\subset L_{0}=\A ^{2k},\]
 where \( L_{i}/L_{i+1}\simeq \wedge ^{2i}E\otimes S^{k-i}TM \), whereas for
\( n=2k+1 \), \( k\geq 0 \), 
\[
L_{k}\subset L_{k-1}\subset \cdots \subset L_{1}\subset L_{0}=\A ^{2k+1},\]
 where \( L_{i}/L_{i+1}\simeq \wedge ^{2i+1}E\otimes S^{k-i}TM_{0} \) (in both
cases we set \( L_{k+1}=0 \)). Thus, \( \A ^{1}=E \), \( \A ^{2}={\mathbb {A}} \)
fits in the Atiyah sequence (\ref{eqn:atiyah}), \( \A ^{3} \) fits in the
exact sequence 
\[
0\rightarrow \wedge ^{3}E\rightarrow \A ^{3}\rightarrow E\otimes TM_{0}\rightarrow 0,\]
 and so on. Such a description in terms of quotients and filtered bundles is
clearly unwieldy to use in practice, so one usually has to resort to local coordinates.
Another alternative is to split the Atiyah sequence (\ref{eqn:atiyah}) by fixing
a linear connection \( \nabla  \) on \( E \), preserving \( \bform  \). Then
\( \A ^{\cdot } \) can be identified with \( \wedge ^{\cdot }E\otimes S^{\cdot }TM_{0} \),
with the nonzero Poisson brackets of generators \( f\in C^{\infty }(M_{0}) \)
(degree \( 0 \)), \( X,Y\in \Gamma (E) \) (degree \( 1 \)) and \( v,w\in \Gamma (TM) \)
(degree \( 2 \)) given as follows: 
\[
\begin{array}{cclcccl}
\{v,f\} & = & v\cdot f &  & \{v,X\} & = & \nabla _{v}X\\
\{X,Y\} & = & <X,Y> &  & \{v,w\} & = & [v,w]+\omega (v,w)
\end{array}\]
 where \( [\cdot ,\cdot ] \) denotes the Lie bracket of vector fields and \( \omega \in \Omega ^{2}(M_{0},\wedge ^{2}E) \)
is the curvature of \( \nabla  \). The Jacobi identity for \( \{\cdot ,\cdot \} \)
is then a consequence of the fact that \( \nabla  \) preserves \( \bform  \),
and of the Bianchi identity.

\section{\label{sec:CA} Poisson manifolds and Courant algebroids.}

Let us now turn our attention to NQ-structures on symplectic N-manifolds, i.e.
homological vector fields of degree \( +1 \) preserving the symplectic form.
By Lemma \ref{lemma:exact}, such a vector field is necessarily Hamiltonian.
Thus, if our manifold \( M \) is of degree \( n \), we can speak of \emph{Hamiltonian
\( n \)-algebroids}.

Consider the case \( n=1 \) first. By Proposition \ref{prop:deg1} \( M=M_{1} \)
is of the form \( T^{*}[1]M_{0} \) for some ordinary manifold \( M_{0} \),
with \( \Omega =d\theta _{i}dx^{i} \), and the corresponding odd Poisson bracket
is just the Schouten bracket \( [\cdot ,\cdot ] \) of multivector fields. By
the above, any NQ-structure on \( M \) is determined by a quadratic Hamiltonian
\( \pi =-\frac{1}{2}\pi ^{ij}(x)\theta _{i}\theta _{j} \) satisfying \( [\pi ,\pi ]=0 \),
i.e. a Poisson bivector field on \( M_{0} \). Thus, 

\begin{prop}
Symplectic NQ-manifolds of degree \( 1 \) are in 1-1 correspondence with ordinary
Poisson manifolds.
\end{prop}
The Poisson bracket of functions on \( M_{0} \) can be expressed in terms of
\( \pi  \) and the Schouten bracket as the so-called \emph{derived bracket:
\[
\{f,g\}=[[f,\pi ],g]\]
} By nondegeneracy of \( [\cdot ,\cdot ] \) this completely determines \( \pi  \):
\( \pi ^{ij}=\{x^{i},x^{j}\} \). The homological Hamiltonian vector field \( d_{\pi }=[\pi ,\cdot ] \)
is the differential in the standard complex computing the Poisson cohomology
of \( M_{0} \). 

Consider now the case \( n=2 \). By Theorem \ref{thm:deg2}, such a manifold
\( (M,\Omega ) \) corresponds to a pseudo-Euclidean vector bundle \( E \)
over some ordinary manifold \( M_{0} \). By Lemma \ref{lemma:exact}, an NQ-structure
on \( M \) is determined by a \emph{cubic} Hamiltonian \( \Theta  \) (``BRST
charge'') satisfying the \emph{structure equation} 
\begin{equation}
\label{eqn:structure}
\{\Theta ,\Theta \}=0
\end{equation}
 where \( \{\cdot ,\cdot \} \) is the even Poisson bracket on \( M \) of degree
\( -2 \). One should expect such a \( \Theta  \) to correspond to some structure
on \( E \), recoverable via derived brackets. This structure is called a \emph{Courant
algebroid}. It was first defined in \cite{LWX1}. We give an equivalent definition
from \cite{Roy1}:

\begin{definition}
A \emph{Courant algebroid} is a pseudo-Euclidean vector bundle \\
\( (E,\bform ) \) over a manifold \( M_{0} \), together with a bilinear operation
\( \circ  \) on \( \Gamma (E) \) and a bundle map \( a:E\rightarrow TM_{0} \)
(the anchor), satisfying the following properties:
\end{definition}
\begin{enumerate}
\item \( e\circ (e_{1}\circ e_{2})=(e\circ e_{1})\circ e_{2}+e_{1}\circ (e\circ e_{2})\; \; \forall e,e_{1},e_{2}\in \Gamma (E) \);
\item \( a(e_{1}\circ e_{2})=[a(e_{1}),a(e_{2})]\; \; \forall e_{1},e_{2}\in \Gamma (E) \);
\item \( e_{1}\circ (fe_{2})=f(e_{1}\circ e_{2})+(a(e_{1})\cdot f)e_{2}\; \; \forall e_{1},e_{2}\in \Gamma (E),\; f\in C^{\infty }(M_{0}) \);
\item \( <e,e_{1}\circ e_{2}+e_{2}\circ e_{1}>=a(e)\cdot <e_{1},e_{2}>\; \; \forall e,e_{1},e_{2}\in \Gamma (E) \);
\item \( a(e)\cdot <e_{1},e_{2}>=<e\circ e_{1},e_{2}>+<e_{1},e\circ e_{2}>\; \; \forall e,e_{1},e_{2}\in \Gamma (E) \).
\end{enumerate}
One can think of sections of \( E \) as acting on \( E \) by ``left multilications''
\( e\circ  \). Properties (1), (3) and (5) mean that this action preserves
the whole structure, while property (2) says that the anchor is equivariant
with respect to this action and the Lie derivative action on \( TM_{0} \).
Property (1) defines a Leibniz algebra structure on \( \Gamma (E) \); if \( \circ  \)
were skew-symmetric, it would be a Lie algebra. Property (4) says that the symmetric
part of \( \circ  \) is ``infinitesimal''. If we define the \emph{co-anchor}
\( a^{*}:T^{*}M_{0}\rightarrow E \) by \( <a^{*}\nu ,e>=\nu (a(e)) \), and
set \( D=a^{*}d \), property (4) reads:
\[
e_{1}\circ e_{2}+e_{2}\circ e_{1}=D<e_{1},e_{2}>\]
 Together with property (2), this implies the useful identity 
\begin{equation}
\label{eqn:aastar}
aa^{*}=0
\end{equation}
 implying that the image of \( a^{*} \) is isotropic. Other useful identities
implied by the definition are 
\begin{equation}
\label{eqn:action_exact}
\begin{array}{rcl}
e\circ Df & = & D<e,Df>\\
Df\circ e & = & 0
\end{array}
\end{equation}
 and more generally 
\begin{equation}
\label{eqn:action_all}
\begin{array}{rcl}
e\circ a^{*}\nu  & = & a^{*}L_{a(e)}\nu \\
a^{*}\nu \circ e & = & -a^{*}\iota _{a(e)}d\nu 
\end{array}
\end{equation}

\begin{example}
Let \( M_{0}=\textrm{point} \). Then \( E \) is just a pseudo-Euclidean vector
space, and a Courant algebroid structure is just a Lie algebra such that \( \bform  \)
is ad-invariant. More generally, any Courant algebroid with \( a=0 \) is a
bundle of such Lie algebras. By (\ref{eqn:aastar}), this must always be the
case if \( \bform  \) is definite.
\end{example}

\begin{example}
\label{example:standardCA} Let \( E=TM_{0}\oplus T^{*}M_{0} \) with the canonical
\( \bform  \) as in Example \ref{example:T^*[2]}. Define \( a \) to be the
projection onto \( TM_{0} \), and \( \circ  \) by 
\begin{equation}
\label{eqn:standard_brack}
(v,\xi )\circ (w,\eta )=([v,w],L_{v}\eta -\iota _{w}d\xi )
\end{equation}
 This operation (rather, its skew-symmetrization) was used originally by T.
Courant \cite{Cou} to describe integrable Dirac structures on \( M_{0} \),
including closed 2-forms, Poisson bivectors and foliations. This example is
responsible for the term ``Courant algebroid'', coined by the authors of \cite{LWX1}
who considered generalizations of this example to dual pairs of Lie algebroids.
\end{example}
\begin{theorem}
\label{thm:CA} Symplectic NQ-manifolds of degree \( 2 \) are in 1-1 correspondence
with Courant algebroids. 
\end{theorem}
\begin{proof}
Let \( (M,\Omega ) \) be the symplectic N-manifold corresponding to \\ \( (E,\bform ) \),
with \( (\A ^{\cdot },\{\cdot ,\cdot \}) \) its graded Poisson algebra of polynomial
functions. Then \( \A ^{0} \) is identified with \( C^{\infty }(M_{0}) \),
\( \A ^{1} \) -- with \( \Gamma (E) \), and \( \{\cdot ,\cdot \} \) restricted
to \( \A ^{1} \) is just \( \bform  \). Let \( \Theta \in \A ^{3} \) satisfy
the structure equation (\ref{eqn:structure}). Given arbitrary \( f\in \A ^{0} \)
and \( e,e_{1},e_{2}\in \A ^{1} \), define \( a \) and \( \circ  \) as the
derived brackets 
\begin{equation}
\label{eqn:derived_brackets}
\begin{array}{rcl}
a(e)\cdot f & = & \{\{e,\Theta \},f\}=\{\{\Theta ,f\},e\}\\
e_{1}\circ e_{2} & = & \{\{e_{1},\Theta \},e_{2}\}
\end{array}
\end{equation}
 The proof that this defines a Courant algebroid is exactly the same as the
proof of Theorem 3.7.3 in \cite{Roy1}, which asserted the same statement for
\( E=A\oplus A^{*} \), \( M=T^{*}[2]A[1] \). But the splitting plays no role
in the argument, so we will not repeat it here. It is a straightforward derivation
using the basic properties of \( \{\cdot ,\cdot \} \) and the structure equation. 

Due to the nondegeneracy of \( \{\cdot ,\cdot \} \), the above derived brackets
uniquely determine \( \Theta  \). In an affine Darboux chart \( (q^{i},\xi ^{a},p_{i}) \)
on \( M \), corresponding to a chart \( \{x^{i}\} \) on \( M_{0} \) and a
local basis \( \{e_{a}\} \) of sections of \( E \) such that \( <e_{a},e_{b}>=g_{ab}=\textrm{const}. \),
\( \Omega =dp_{i}dq^{i}+\frac{1}{2}d\xi ^{a}g_{ab}d\xi ^{b} \), and \( \Theta \in \A ^{3} \)
is of the form 
\begin{equation}
\label{eqn:CourantTensor}
\Theta =\xi ^{a}A_{a}^{i}(q)p_{i}-\frac{1}{6}\phi _{abc}(q)\xi ^{a}\xi ^{b}\xi ^{c}
\end{equation}
 It clearly follows from (\ref{eqn:derived_brackets}) that 
\begin{equation}
\label{eqn:components}
\begin{array}{rcl}
A^{i}_{a} & = & a(e_{a})\cdot x^{i}\\
\phi _{abc} & = & <e_{a}\circ e_{b},e_{c}>
\end{array}
\end{equation}
 where \( e_{a}=g_{ab}\xi ^{b} \). So, conversely, given a Courant algebroid
structure on \\
\( (E,\bform ) \), we can define \( \Theta  \) in a Darboux chart by the formulas
(\ref{eqn:CourantTensor}), (\ref{eqn:components}). It is an easy check, using
the coordinate transformations (\ref{eqn:coord_change}), that \( \Theta  \)
is in fact globally defined and obeys (\ref{eqn:structure}) due to the properties
of a Courant algebroid. In fact, 
\[
\begin{array}{rcl}
\{\Theta ,\Theta \} & = & <a^{*}dx^{i},a^{*}dx^{j}>p_{i}p_{j}+\xi ^{a}\xi ^{b}dx^{j}([a(e_{a}),a(e_{b})]-a(e_{a}\circ e_{b}))p_{j}+\\
 & + & \frac{1}{12}<(e_{a}\circ e_{b})\circ e_{c}+e_{b}\circ (e_{a}\circ e_{c})-e_{a}\circ (e_{b}\circ e_{c}),e_{d}>\xi ^{a}\xi ^{b}\xi ^{c}\xi ^{d}
\end{array}\]
 So \( \{\Theta ,\Theta \}=0 \) if and only if \( (E,\bform ,a,\circ ) \)
is a Courant algebroid.
\end{proof}
\begin{example}
If \( M_{0} \) is a point, \( \Theta \in \wedge ^{3}E^{*} \) is the Cartan
structure tensor corresponding to a Lie algebra structure on \( E \) making
\( \bform  \) ad-invariant: \\
\( \Theta (X,Y,Z)=<[X,Y],Z> \).
\end{example}

\begin{example}
For \( E=TM_{0}\oplus T^{*}M_{0} \) with the standard Courant algebroid structure,
\( M=T^{*}[2]T[1]M_{0} \), and \( \Theta =\Theta _{0} \) is the Hamiltonian
corresponding to the de Rham vector field \( d \) on \( T[1]M_{0} \). This
example will be studied in detail in the next section.
\end{example}

\begin{example}
More generally, for \( E=A\oplus A^{*} \), the double of a Lie bialgebroid
\( (A,A^{*}) \) \cite{LWX1}, \( M=T^{*}[2]A[1]\simeq T^{*}[2]A^{*}[1] \),
and \( \Theta =\mu +\gamma  \) where \( \mu  \) and \( \gamma  \) are the
commuting Hamiltonians defining the Lie bialgebroid structure \cite{Roy1}.
\end{example}

\begin{example}
In particular, let a Lie algebra \( {\mathfrak {g}} \) act on a manifold \( M_{0} \)
(on the right). Endow the trivial bundle \( M\times {\mathfrak {g}} \) with
the transformation (semidirect product) Lie algebroid structure, the dual bundle
with the zero structure, and let \( E=M\times ({\mathfrak {g}}\oplus {\mathfrak {g}}^{*}) \)
be the double. Then \( M=({\mathfrak {g}}\oplus {\mathfrak {g}}^{*})[1]\times T^{*}[2]M_{0} \),
and the BRST charge is the classical one: 
\[
\Theta =\xi ^{a}v_{a}^{i}(q)p_{i}-\frac{1}{2}\xi ^{a}\xi ^{b}C_{ab}^{c}\theta _{c}\]
 where \( v_{a}=v_{a}^{i}(x)\frac{\partial }{\partial x^{i}} \) are the infinitesimal
generators of the action corresponding to a basis \( (e_{a}) \) of \( {\mathfrak {g}} \),
and \( (C_{ab}^{c}) \) are the structure constants in this basis. Notice that
our grading is different from the ``ghost degree'' grading. This leads to
a different grading in the cohomology as well, only the parity remains the same.
However, under sufficient regularity assumption, the total cohomology is isomorphic
to the algebra of functions on the reduced phase space, as proved in \cite{KoSt},
so this change in grading is of no consequence.
\end{example}
Given a Courant algebroid, \( (\A ^{\cdot },D=\{\Theta ,\cdot \}) \) becomes
a differential complex, the \emph{standard complex} of \( (E,\bform ,a,\circ ) \).
Since \( D \) is a Hamiltonian vector field and \( M \) has degree \( 2 \),
it is completely determined by its action on \( \A ^{0}=C^{\infty }(M_{0}) \)
and \( \A ^{1}=\Gamma (E) \). In fact, the formulas (\ref{eqn:derived_brackets})
imply that \( D:C^{\infty }(M_{0})\rightarrow \Gamma (E) \) is just \( D=a^{*}d \),
while \( D:\Gamma (E)\rightarrow \Gamma ({\mathbb {A}}) \) is given by 
\[
(De_{1})\cdot e_{2}=e_{1}\circ e_{2}\]
 where \( \cdot  \) denotes the Lie algebroid action of \( \Gamma ({\mathbb {A}}) \)
on \( \Gamma (E) \). One interesting feature of this complex is that it is
infinite in general.

The cohomology groups \( H^{\cdot }(M,D) \) in low degrees have familiar structural
interpretations. Thus, \( H^{0} \) is the space of smooth functions on \( M_{0} \)
that are constant along the leaves of the anchor foliation; it is equal to \( {\mathbb {R}} \)
for transitive Courant algebroids (those with surjective anchor). \( H^{1} \)
is the space of sections of \( E \) acting trivially on \( E \), modulo those
of the form \( Df \) for some function \( f \). Further, \( H^{2} \) is the
space of linear vector fields on \( E \) preserving the Courant algebroid structure,
modulo those generated by sections of \( E \) as \( e\circ  \). \( H^{3} \)
is the space of infinitesimal deformations of the Courant algebroid structure,
modulo the trivial ones generated by \( \Gamma ({\mathbb {A}}) \), while \( H^{4} \)
houses the obstructions to extending an infinitesimal deformation to a formal
one. For the purposes of deformation theory it is customary to shift the grading
by \( -2 \), the degree of \( \{\cdot ,\cdot \} \), so that \( H^{3} \) and
\( H^{4} \) become, respectively, \( H^{1} \) and \( H^{2} \). Notice the
analogy with Poisson manifolds, and with Lie algebras.

\section{\label{sec:exact} Higher de Rham complexes and exact Courant algebroids .}

Let us consider the case \( E=TM_{0}\oplus T^{*}M_{0} \). In this case \( M=T^{*}[2]T[1]M_{0}=T^{*}\Pi TM_{0} \),
as in Example \ref{example:T^*[2]}. Any local coordinates \( (x^{i}) \) on
\( M_{0} \) give rise to affine Darboux coordinates \( (q^{i},\xi ^{i},p_{i},\theta _{i}) \),
so that \( \Omega =dq^{i}dp_{i}+d\xi ^{i}d\theta _{i} \). This manifold has
a rich additional structure. First of all, there is a canonical isomorphism,
the \emph{Legendre transformation}, that identifies \( T^{*}\Pi TM_{0} \) with
\( T^{*}\Pi T^{*}M_{0} \) by exchanging the momenta conjugate to the fibre
coordinates on \( \Pi TM_{0} \) with the fibre coordinates on \( \Pi T^{*}M_{0} \),
and vice versa. This holds for any vector bundle \( A\rightarrow M_{0} \) \cite{Roy1}.
Thus, \( M \) is what is known as a \emph{double vector bundle}: 
\[
\begin{array}{ccc}
M & \longrightarrow  & T^{*}[1]M_{0}\\
\downarrow  &  & \downarrow \\
T[1]M_{0} & \longrightarrow  & M_{0}
\end{array}\]
 This means that each arrow is a vector bundle; the horizontal pair of arrows
is a morphism of the vertical vector bundles, while the vertical pair is a morphism
of the horizontal bundles. Because of this structure, the algebra \( \A ^{\cdot } \)
of polynomial functions on \( M \) acquires a double grading; we shall use
the notation \( \A ^{\cdot ,\cdot }=\oplus _{k,l\geq 0}\A ^{k,l} \). The Euler
vector fields corresponding to the two gradings are, respectively, 
\[
\epsilon _{1}=p_{i}\frac{\partial }{\partial p_{i}}+\theta _{i}\frac{\partial }{\partial \theta _{i}}\]
 and 
\[
\epsilon _{2}=p_{i}\frac{\partial }{\partial p_{i}}+\xi ^{i}\frac{\partial }{\partial \xi ^{i}}\]
 The N-manifold structure on \( M=T^{*}[2]T[1]M_{0}=T^{*}[2]T^{*}[1]M_{0} \)
corresponds to the total grading 
\[
\epsilon =\epsilon _{1}+\epsilon _{2}=2p_{i}\frac{\partial }{\partial p_{i}}+\xi ^{i}\frac{\partial }{\partial \xi ^{i}}+\theta _{i}\frac{\partial }{\partial \theta _{i}}\]
 Evidently, the canonical sympelctic form \( \Omega  \) has bi-degree \( (1,1) \),
the corresponding Poisson bracket -- \( (-1,-1) \). 

There is another special canonical isomorphism: \( T^{*}[2]T^{*}[1]M_{0}=T^{*}\Pi T^{*}M_{0}\simeq \Pi T\Pi T^{*}M_{0}=T[1]T^{*}[1]M_{0} \).
It comes from treating the conjugate momenta \( p_{i},\xi ^{i} \) as the differentials
of the basic coordinates: \( p_{i}=d\theta _{i} \), \( \xi ^{i}=dx^{i} \).\footnote{
For a thorough treatment of double vector bundle and canonical isomorphisms
(without parity shifts), see \cite{KonUrb}
} Correspondingly, on \( M \) there is the canonical (de Rham) vector field
\[
D=\xi ^{i}\frac{\partial }{\partial q^{i}}+p_{i}\frac{\partial }{\partial \theta _{i}}\]
 This vector field is homological and Hamiltonian: \( D=\{\Theta _{0},\cdot \} \),
where \( \Theta _{0}=\xi ^{i}p_{i}\in \A ^{1,2} \) obeys (\ref{eqn:structure}).
With respect to the fibration \( M=T^{*}\Pi TM_{0}\rightarrow \Pi TM_{0} \),
it is just the Hamiltonian lift of the de Rham vector field \( d=\xi ^{i}\frac{\partial }{\partial x^{i}} \)
on \( T[1]M_{0} \). With respect to the other fibration \( M=T^{*}\Pi T^{*}M_{0}\rightarrow \Pi T^{*}M_{0} \),
\( \Theta _{0} \) is quadratic and generates the canonical Schouten bracket
on \( T^{*}[1]M_{0} \) as a derived bracket: 
\begin{equation}
\label{eqn:schouten}
[X,Y]=\{\{X,\Theta _{0}\},Y\}
\end{equation}
 for multivector fields \( X,Y \). 

The vector field \( D \) has bi-degree \( (0,+1) \). This can be expressed
by the relations 
\[
\begin{array}{rcl}
{[\epsilon _{1},D]} & = & 0\\
{[\epsilon _{2},D]} & = & D
\end{array}\]
 This means that \( D \) preserves the subspace \( \A ^{k,\cdot } \) for each
\( k\geq 0 \), and the restriction to each \( \A ^{k,\cdot } \) is a differential.
For \( k=0 \), \( (\A ^{0,\cdot },D)=(\Omega ^{\cdot }(M_{0}),d) \), the de
Rham complex of \( M_{0} \); we call the complexes \( (\A ^{k,\cdot },D) \)
for \( k>0 \) the \emph{higher de Rham complexes}. In particular, \( \A ^{1,\cdot }[1] \)
is naturally identified with vector fields on \( T[1]M_{0} \), i.e. derivations
of \( \Omega ^{\cdot }(M_{0}) \), with \( D=[d,\cdot ] \). In general, \( \A ^{k,l} \)
is the space of symbols of differential operators on \( T[1]M_{0} \) of order
\( k \) and degree \( l-k \).

What is important is that \( \epsilon _{1} \), viewed as a cocycle with respect
to \( [D,\cdot ] \), is actually a coboundary. That is, there exists a vector
field \( \iota  \) on \( M \), of bi-degree \( (0,-1) \), such that 
\begin{equation}
\label{eqn:homotopy}
[D,\iota ]=D\iota +\iota D=\epsilon _{1}
\end{equation}
 In coordinates, 
\[
\iota =\theta _{i}\frac{\partial }{\partial p_{i}}\]
 It is clear that \( \iota  \) satisfies the above relation, and is homological.
With respect to the symplectic structure \( \Omega  \) on \( M \), \( \iota  \)
is not Hamiltonian, but rather satisfies the relation: 
\[
L_{\iota }\Omega +\omega =0,\]
 where \( \omega =d\theta _{i}dx^{i} \) is the canonical odd symplectic form
on \( T^{*}[1]M_{0} \), pulled back to \( M \). This relates symplectic N-manifolds
of degrees 1 and 2.

To show that \( \iota  \) is well-defined, one can check that it is invariant
under coordinate transformations, but we shall provide a coordinate-free description
instead. It suffices to describe the action of \( \iota  \) on \( \A ^{0,\cdot } \)
and \( \A ^{1,\cdot } \). On \( \A ^{0,\cdot } \) it acts trivially; let \( V\in \A ^{1,l} \).
Then \( V \) is a derivation of \( \Omega ^{\cdot }(M_{0}) \) of degree \( l-1 \).
The restriction of \( V \) to \( \Omega ^{0}(M_{0})=C^{\infty }(M_{0}) \)
is a derivation with values in \( \Omega ^{l-1}(M_{0}) \), i.e. a \( (l-1) \)-form-valued
vector field \( \bar{V} \) on \( M_{0} \). We set \( \iota (V)=\iota _{\bar{V}} \),
the contraction operator. Now, \( \epsilon _{1} \) acts on \( \A ^{1,\cdot } \)
as the identity operator. Therefore, by (\ref{eqn:homotopy}), the operators
\( \iota D \) and \( D\iota  \) are complementary projections. This leads
to the following decomposition: 
\begin{equation}
\label{eqn:decompos}
\A ^{1,\cdot }[1]=\A ^{1,\cdot +1}=\textrm{Der}^{\cdot }(\Omega (M_{0}))\simeq {\mathcal{G}}[\tau ]={\mathcal{G}}^{\cdot }[-1]\oplus {\mathcal{G}}^{\cdot }
\end{equation}
 as graded vector spaces. Here \( {\mathcal{G}}^{\cdot }=\Gamma (\wedge ^{\cdot }T^{*}M_{0}\otimes TM_{0}) \),
the space of vector-valued differential forms, and \( {\mathcal{G}}[\tau ] \)
denotes the adjunction of a formal odd variable of degree \( -1 \); the differential
\( D=[d,\cdot ] \) acts as \( \frac{d}{d\tau } \). The space \( {\mathcal{G}}^{\cdot } \)
appears in the first summand as ``contractions'', in the second -- as ``Lie
derivatives''. The above decomposition (\ref{eqn:decompos}) is a classic result
of Fr\"{o}licher and Nijenhuis \cite{FroNij}. By viewing \( D \) and \( \iota  \)
as vector fields on \( M \), we obtain the following generalization: 

\begin{prop}
\label{prop:acyclicity} The higher de Rham complexes \( (\A ^{k,\cdot },D) \),
\( k>0 \), are acyclic.
\end{prop}
\begin{proof}
Immediate from (\ref{eqn:homotopy}), since \( \epsilon _{1} \) acts on \( \A ^{k,\cdot } \)as
multiplication by \( k \). This also gives rise to a decomposition of \( \A ^{k,\cdot } \)
generalizing (\ref{eqn:decompos}). 
\end{proof}
The above proposition has several consequences. For example, in terms of the
standard Lie algebroid structure on \( TM_{0} \), the de Rham complex \( (\A ^{0,\cdot },d) \)
is the standard complex with coefficient in the trivial module, whereas the
higher complex \( (\A ^{k,\cdot },D) \) should be viewed as ``the standard
complex with coefficients in the \( k \)the exterior power of the adjoint module'',
by analogy with Lie algebras. Thus, for instance, the complex \( (\A ^{1,\cdot }[1],[d,\cdot ]) \)
is responsible for the deformation theory of this algebroid, and its acyclicity
implies the rigidity of the standard Lie algebroid:

\begin{cor}
\label{cor:rigidityTM} The standard Lie algebroid structure on \( TM_{0} \)
is rigid. 
\end{cor}
The complex \( (\A ^{2,\cdot },D) \) is also interesting. As shown in \cite{Roy1},
1-cocycles in that complex are responsible for Lie algebroid structures on \( T^{*}M_{0} \),
compatible with the standard one on \( TM_{0} \) in the sense that \( (TM_{0},T^{*}M_{0}) \)
is a Lie bialgebroid \cite{MacXu}. More precisely, a \( \gamma \in \A ^{2,1} \)
gives rise to a Lie algebroid on \( T^{*}M_{0} \) if and only if \( \{\gamma ,\gamma \}=0 \),
and it is compatible with \( TM_{0} \) if and only if \( D\gamma =\{\Theta _{0},\gamma \}=0 \).
By Proposition \ref{prop:acyclicity}, \( \gamma =D\pi =\{\Theta _{0},\pi \} \),
where \( \pi \in \A ^{2,0}=\Gamma (\wedge ^{2}TM_{0}) \) is some bivector field.
Then by (\ref{eqn:schouten}), \( \{\gamma ,\gamma \}=0 \) is equivalent to
\( \{\Theta _{0},[\pi ,\pi ]\}=0 \); by Proposition \ref{prop:acyclicity}
(with \( k=3 \)), this can only happen if \( [\pi ,\pi ]=0 \). Therefore,
we have another 

\begin{cor}
\label{cor:Poisson} Any Lie algebroid structure on \( T^{*}M_{0} \) compatible
with the standard one on \( TM_{0} \) comes from a Poisson structure on \( M_{0} \).
\end{cor}
Finally, let us consider \( (M,\Omega ,\Theta _{0}) \) as a symplectic NQ-manifold.
By Theorem \ref{thm:CA}, it must corresond to a Courant algebroid structure
on \( E=TM_{0}\oplus T^{*}M_{0} \); it is easy to see that this Courant algebroid
is the standard one, as in Example \ref{example:standardCA}. Let us denote
this structure by \( E_{0} \). The standard complex of \( E_{0} \) is the
total complex \( (\A ^{\cdot },D) \). The Proposition \ref{prop:acyclicity}
immediately implies

\begin{cor}
\label{cor:standardCA_cohomology} The cohomology of the standard complex of
the Courant algebroid \( E_{0} \) is isomorphic to the de Rham cohomology of
\( M_{0} \).
\end{cor}
We can describe the action of the de Rham cohomology explicitly. Thus, 1-forms
on \( M_{0} \) act on any Courant algebroid over \( M_{0} \) via (\ref{eqn:action_all}),
and closed 1-forms act trivially; in case of \( E_{0} \), they are the only
sections that do. Furthermore, a 2-form \( \beta \in \A ^{0,2}\subset \A ^{2} \)
acts by changing the isotropic splitting of \( E=TM_{0}\oplus T^{*}M_{0} \)
via its graph: \( (v,\xi )\mapsto (v,\xi +\iota _{v}\beta ) \). This does not
change \( \bform  \) and \( a \), but \( \circ  \) in (\ref{eqn:standard_brack})
is modified by the term \( d\beta (v,w,\cdot ) \). So the Courant algebroid
structure is preserved if \( d\beta =0 \); if \( \beta =d\alpha  \), the action
is \( \alpha \circ  \), as in (\ref{eqn:action_all}).

Lastly, \( H^{3}(M_{0},{\mathbb {R}}) \) corresponds to nontrivial deformations
of \( E_{0} \). A 3-form \( \phi  \) acts by adding the term \( \phi (v,w,\cdot ) \)
to the expression (\ref{eqn:standard_brack}) for \( \circ  \). The corresponding
structure is 
\[
\Theta _{\phi }=\Theta _{0}-\phi =\xi ^{i}p_{i}-\frac{1}{6}\phi _{ijk}(q)\xi ^{i}\xi ^{j}\xi ^{k}\]
 where \( \phi =\frac{1}{6}\phi _{ijk}(x)dx^{i}\wedge dx^{j}\wedge dx^{k} \).
It is easy to check that \( \Theta _{\phi } \) obeys the structure equation
(\ref{eqn:structure}) if and only if \( d\phi =0 \). Let us denote the corresponding
Courant algebroid \( E_{\phi } \). Cohomologous 3-forms give rise to isomorphic
Courant algebroids: if \( \phi '-\phi =d\beta  \), changing the splitting by
\( \beta  \) maps \( E_{\phi } \) to \( E_{\phi '} \). In particular, exact
3-forms yield trivial deformations of \( E_{0} \). Corollary \ref{cor:standardCA_cohomology}
implies that these \( E_{\phi } \)'s are the only possible deformations of
\( E_{0} \). 

In effect, we have recovered, via deformation theory, \v{S}evera's classification
of so-called \emph{exact Courant algebroids} \cite{Sev}. These are Courant
algebroids such that the sequence 
\[
0\rightarrow T^{*}M_{0}\stackrel{a^{*}}{\rightarrow }E\stackrel{a}{\rightarrow }TM_{0}\rightarrow 0\]
 is exact. Choosing and isotropic splitting \( \sigma :TM_{0}\rightarrow E \)
(a ``connection'') identifies \( E \) with \( E_{\phi } \), where \( \phi  \)
is the ``curvature'' of \( \sigma  \), given by \( \phi (u,v,w)=<\sigma (u)\circ \sigma (v),\sigma (w)> \).
The cohomology class of \( \phi  \) is a characteristic class, the \emph{\v{S}evera
class} of \( E \). 

\begin{example}
(due to Anton Alexeev) As an example of an exact Courant algebroid with a nontrivial
\v{S}evera class, consider a Lie group \( G \) with a bi-invariant metric \( K \).
Let \( E=({\mathfrak {g}}\oplus {\mathfrak {g}})\times G \), where \( {\mathfrak {g}}=\textrm{Lie}(G) \).
Define \( \bform  \) by 
\[
<(X,Y),(X',Y')>=K(X,X')-K(Y,Y')\]
 and an anchor \( a:E\rightarrow TG \) by 
\[
a(X,Y)=X^{l}-Y^{r}\]
 i.e. \( a(X,Y,g)=gX-Yg \). Finally, define \( \circ  \) on constant sections
by 
\[
(X,Y)\circ (X',Y')=([X,X'],[Y,Y'])\]
 By trivializing both \( TG \) and \( T^{*}G \) as \( G\times {\mathfrak {g}} \)
by left translation and using \( K \), it is easy to see that this defines
an exact Courant algebroid over \( G \). Let us split \( a \) by \( \sigma (Z,g)=\frac{1}{2}(Z,-\textrm{Ad}_{g}Z) \).
One computes the curvature to be \( \phi (X,Y,Z)=-\frac{1}{2}K([X,Y],Z) \),
the Cartan structure tensor. Thus, the \v{S}evera class of this Courant algebroid
is the canonical class of \( G \).
\end{example}


\end{document}